\documentclass[12pt]{amsart}

\usepackage{amsthm}
\usepackage{amsfonts}

\newtheorem{theorem}{Theorem}[section]
\newtheorem{corollary}[theorem]{Corollary}
\newtheorem{lemma}[theorem]{Lemma}
\newtheorem{proposition}[theorem]{Proposition}

\newtheorem{definition}[theorem]{Definition}

\newcommand{\C}{{\mathbb{C}}}

\newcommand{\R}{{\mathbb{R}}}
\newcommand{\Z}{{\mathbb{Z}}}
\newcommand{\Q}{{\mathbb{Q}}}

\numberwithin{equation}{section}
\begin{document}

\title[Equivariant structure constants for weighted projective space]{Equivariant structure constants for ordinary and weighted projective space}
\author{Julianna S. Tymoczko}
\address{Department of Mathematics, University of Iowa, Iowa City, 
IA 52242}
\email{tymoczko@math.uiowa.edu}

\thanks{The author is grateful to Tony Bahri, Linda Chen, Charlie Frohman, Bill Fulton, and Tara Holm for useful conversations.}


\begin{abstract}
We compute the integral torus-equivariant cohomology ring for weighted projective space for two different torus actions by embedding the cohomology in a sum of polynomial rings $\oplus_{i=0}^n \Z[t_1, t_2,\ldots, t_n]$.  One torus action gives a result complementing that of Bahri, Franz, and Ray.  For the other torus action, each basis class for weighted projective space is a multiple of the basis class for ordinary projective space; we identify each multiple explicitly.  We also give a simple formula for the structure constants of the equivariant cohomology ring of ordinary projective space in terms of the basis of Schubert classes, as a sequence of divided difference operators applied to a specific polynomial.
\end{abstract}

\maketitle

\section{Introduction}

Ordinary projective space is the collection of lines through the origin in $\C^n$; weighted projective space is the collection of curves through the origin in $\C^n$ with fixed trajectories.  More formally, given $\lambda \in \Z_{>0}^{n+1}$ the weighted projective space $w^{\lambda} \mathbb{P}^n$ is the quotient 
\[w^{\lambda} \mathbb{P}^n =  \{[z_0, \ldots, z_n] \in (\C^*)^{n+1}\}/  \sim\}\]
by the equivalence relation $[z_0, z_1, \ldots, z_n] \sim [t^{\lambda_0}z_0, t^{\lambda_1}z_1, \ldots, t^{\lambda_n}z_n]$ for each $t \in \C^*$.

Mathematical physicists and algebraic geometers study weighted projective space because it is the simplest example of an orbifold.  Moreover as the weight $\lambda$ varies, weighted projective spaces give local models for all of the orbifold singularities that correspond to the quotient by a finite {\em abelian} group.  

This paper studies the equivariant cohomology of weighted and ordinary projective space with integer coefficients under several torus actions.  It was inspired by recent work of Bahri, Franz, and Ray in \cite{BFR}, which computes the equivariant cohomology of weighted projective space with respect to the action of the torus $T = (\C^*)^{n+1}$ defined by $(t_0,\ldots,t_n) \cdot [z_0,\ldots, z_n] = [t_0z_0,\ldots, t_nz_n]$.  \footnote{They actually use the quotient $T/(\C^*)^{\lambda}$ but the extra factor stabilizes $w^{\lambda} \mathbb{P}^n$.  Hence the $T$-equivariant cohomology is an algebra over a polynomial ring in $n+1$ variables $y_i$; imposing the relation $\sum_{i=0}^n \lambda_iy_i = 0$ gives the $T/(\C^*)^{\lambda}$-equivariant cohomology.} 

We study a family of torus actions parametrized by $\mu \in \Z_{>0}^{n+1}$.  The action of $T^{\mu}$ is defined by $(t_0,\ldots,t_n) \cdot [z_0,\ldots, z_n] = [t_0^{\mu_0}z_0,\ldots, t_n^{\mu_n}z_n]$.  We show that the integral equivariant cohomology of each weighted projective space $w^{\lambda} \mathbb{P}^n$ with respect to each $T^{\mu}$ can be realized as a subring of $\oplus_{i=0}^n \Z[t_0, \ldots, t_n]$.   

The proof uses localization techniques, in particular the theory of Goresky-Kottwitz-MacPherson for equivariant cohomology \cite{GKM}.  We find $H^*_{T^{\mu}}(w^{\lambda}\mathbb{P}^n; \Q)$ for all $\mu$ and $\lambda$ using the GKM approach.  Franz and Puppe adapted GKM theory to integral equivariant cohomology in \cite{FP1}.  However their tools require conditions that are satisfied only when $\mu= (1,1,1,\ldots)$, as discussed in Section 2.  (This is the case considered in \cite{BFR}.)

When $\mu = \lambda$ we can nonetheless completely identify the integral equivariant cohomology of weighted projective space.

\medskip

\noindent {\bf Theorem.} {\em The equivariant cohomology ring $H^*_{T^{\lambda}}(w^{\lambda} \mathbb{P}^n)$ has a natural basis of scalar multiples of the ordinary Schubert classes in $H^*_T(\mathbb{P}^n)$.}

\medskip

We give this scalar explicitly in Theorem \ref{weighted basis} and show that it is the scalar found for ordinary cohomology by Kawasaki in \cite{K}.  

This reduces the work to calculations in the equivariant cohomology ring of ordinary projective space.  Our final result computes the structure constants with respect to the basis of Schubert classes in $H^*_T(\mathbb{P}^n)$.  

The formula we obtain is remarkably elegant, yet as far as we can determine appears nowhere in the literature.  Presentations of the ring  $H^*_T(\mathbb{P}^n)$ in terms of generators and relations are well known \cite[Lecture 3]{F}, as are Chevalley-Monk formulas that describe the product of a degree-one class with an arbitrary class in $H^*_T(G/P)$ \cite[Lecture 15]{F}.  Knutson-Tao, Kreiman, and Molev each provide a powerful construction of the equivariant structure constants of Grassmannians \cite{KT, Kr, M}.  However, these constructions are unnecessarily complicated for projective space; for instance, it is not even combinatorially evident that they produce the formula provided here!

Let $p_i$ be the generator of $H^{2i}_T(\mathbb{P}^n)$.  The divided difference operator is the function $\partial_i: \C[t_0,t_1, \ldots, t_n] \rightarrow \C[t_0,t_1,\ldots,t_n]$ defined by 
\[\partial_i p( \ldots, t_i, t_{i+1},\ldots) = \frac{p( \ldots, t_{i+1},t_i, \ldots) - p( \ldots, t_i, t_{i+1},\ldots)}{t_i - t_{i+1}}.\]

\medskip

\noindent {\bf Theorem \ref{ordinary equivariant projective}.} {\em For each $i=0,\ldots,n$ let $p_i$ denote the generator of $H^i_T(\mathbb{P}^n)$.  If $i,j,k$ are integers with $1 \leq i < j < k \leq \max\{i+j, n\}$ then the coefficients in $p_i p_j = \sum c_{ij}^k p_k$ are
\[c_{ij}^k = \partial_{k-1} \partial_{k-2} \cdots \partial_{j+1} \partial_j q\]
and $c_{ij}^j = q$, where $q = (t_0-t_j)(t_1-t_j)(t_2-t_j)(t_3-t_j) \cdots (t_{i-1}-t_j)$.}

\medskip

\section{Background and definitions}

\subsection{Weighted projective space} Let $T$ be the complex algebraic torus $\C^* \times \cdots \times \C^*$ considered as the set of invertible diagonal matrices acting on $\C^{n+1}$.  Fix a weight $\lambda \in \textup{Hom}(\C^*, T)$.  The action of $\lambda$ on $t \in \C^*$ can be written $\lambda \cdot t = (t^{\lambda_0}, t^{\lambda_1}, \ldots, t^{\lambda_n})$ for fixed positive integers $\lambda_i $.  We take $\lambda \in \Z^{n+1}_{>0}$.

The weighted projective space $w^{\lambda} \mathbb{P}^n$ is defined to be the quotient of $\C^{n+1} - \{0\}$ by the action of $\lambda(\C^*)$, namely the equivalence classes
\[ w^{\lambda} \mathbb{P}^n = \{[z_0, \ldots, z_n] \in \C^{n+1}- \{0\}\}/  \sim \}\]
under the equivalence relation defined by $[z_0, z_1, \ldots, z_n] \sim [t^{\lambda_0}z_0, t^{\lambda_1}z_1, \ldots, t^{\lambda_n}z_n]$ for each $t \in \C^*$.  For instance, when $\lambda = (1,1,\ldots, 1)$ then $w^{\lambda} \mathbb{P}^n$ is ordinary projective space.  Other authors suppress $\lambda$ in the notation; we retain it because we compare different weights.  Weighted projective space is a complex projective algebraic variety, and in fact a toric variety \cite{O}.
 
The torus $T$ acts on $w^{\lambda} \mathbb{P}^n$ coordinate-wise in many different ways.  We consider a family of $T$-actions parametrized by $\Z^{n+1}$.  If $\mu \in \Z^{n+1}$ the action $T^{\mu}$ is defined so that for each $(t_0, \ldots, t_n) \in T$ and $[z_0,\ldots z_n] \in w^{\lambda} \mathbb{P}^n$
\[(t_0, t_1, \ldots, t_n) \cdot [z_0, z_1, \ldots, z_n] = [t_0^{\mu_0}z_0, t_1^{\mu_1}z_1, \ldots, t_n^{\mu_n}z_n].\]

The $T^{\mu}$-fixed points and one-dimensional $T^{\mu}$-orbits of weighted projective space are constructed in the same way as for ordinary projective space.  The $T^{\mu}$-fixed points in $w^{\lambda} \mathbb{P}^n$ are the points with exactly one nonzero coordinate.  If $e_i$ is the $i^{th}$ standard basis vector of $\C^{n+1}$ then denote the fixed point $[e_i]$ by $x_i$.  There is exactly one one-dimensional $T^{\mu}$-orbit $O_{ij}$ for each pair $i,j$ with $0 \leq i < j \leq n$:
\[O_{ij} = \{[z_0,z_1,\ldots, z_n] \in w^{\lambda} \mathbb{P}^n: z_k=0 \textup{ if } k \neq i,j\}.\]
The closure $\overline{O_{ij}}$ contains $x_i$ and $x_j$ but no other $T^{\mu}$-fixed point.  In general for each subset $I \subseteq \{0,1,\ldots, n\}$ there is a unique $T^{\mu}$-orbit $O_I$ consisting of the points $[z_0,z_1,\ldots, z_n]$ with coordinate $z_i$ nonzero if and only if $i \in I$.

\begin{proposition} \label{Tmu-CW}
Let $T^{\mu}_{\R}$ denote the real subtorus of $T$ with the action weighted by $\mu$, so $T^{\mu}_{\R} \cong (S^1)^{n+1}$.  For each $\lambda$ and $\mu$ the weighted projective space $w^{\lambda} \mathbb{P}^n$ is a $T^{\mu}_{\R}$-CW complex.  
\end{proposition}

\begin{proof}
Associate a simplex $\Delta_I \subseteq \R^n$ to each subset $I \subseteq \{0,1,\ldots,n\}$ by taking the convex hull of the standard basis element $e_i$ for each nonzero $i \in I$, together with the origin if $0 \in I$.  Let $\textup{Stab}(O_I) \subseteq T_{\R}^{\mu}$ denote the stabilizer of $O_I$.  The equivariant cell $\Delta_I \times T_{\R}/\textup{Stab}(O_I)$ carries an action of $T_{\R}$ by (unweighted) group multiplication on the second factor.  Define attaching maps $f_I: \Delta_I \times T_{\R}/\textup{Stab}(O_I) \rightarrow w^{\lambda} \mathbb{P}^{n}$ by
\[({\bf v}, {\bf g}) \mapsto \left[g_0^{\mu_0} \left(1 - \sum_{j=1}^n v_j \right), g_1^{\mu_1} v_1, \ldots, g_n^{\mu_n} v_n\right].\]
Weighted projective space is a $T^{\mu}_{\R}$-CW complex if each $f_I$ is a $T_{\R}^{\mu}$-equivariant map that is a homeomorphism on the interior of each cell.  Equivariance follows because $f_I({\bf t} \cdot ({\bf v}, {\bf g})) = {\bf t}^{\lambda} \cdot f_I({\bf v}, {\bf g})$.  The image under $f_I$ of the interior of $\Delta_I \times T_{\R}/\textup{Stab}(O_I)$ lies in $O_I$ by construction.  Proving that $f_I$ is a homeomorphism between $O_I$ and the interior of the $I^{th}$ cell is a routine exercise left to the reader.
\end{proof}

\subsection{Equivariant cohomology} Let $X$ be a topological space with a suitable action of a torus $T$.  GKM (Goresky-Kottwitz-MacPherson) theory gives conditions for the inclusion of fixed points $X^T \hookrightarrow X$ to induce an injection $H^*_T(X) \rightarrow H^*_T(X^T)$.  To identify the image explicitly, GKM theory constructs a combinatorial graph from the data of $T$-fixed points and one-dimensional $T$-orbits of $X$, and then uses an algebraic algorithm to compute equivariant cohomology from this graph.  GKM theory uses $\Q$ coefficients \cite{GKM}; in this section we use \cite{FP1} to modify GKM theory for $\Z$ coefficients.  Unless otherwise stated $H^*_T(X)$ means $H^*_T(X; \Z)$.

\begin{lemma} \label{equivariantly formal}
For each $\mu$ and each $\lambda$ the weighted projective space $w^{\lambda} \mathbb{P}^n$ is $T^{\mu}$-equivariantly formal, in the sense that $H^*_{T^{\mu}}(w^{\lambda} \mathbb{P}^n)$ is free over $H^*(BT^{\mu})$.  If $T^{\mu}_{\R}$ denotes the real subtorus of $T^{\mu}$ then $w^{\lambda} \mathbb{P}^n$ is also $T^{\mu}_{\R}$-equivariantly formal.
\end{lemma}

\begin{proof}
The proof is from \cite[Lemma 2.2]{BFR}.  We reproduce it here to show that it is independent of $\mu$ and of the choice of $T^{\mu}$ or $T^{\mu}_{\R}$.  Kawasaki showed that the ordinary cohomology of $w^{\lambda} \mathbb{P}^n$ is free and vanishes in odd degrees.  This means that $E_2 = E_{\infty}$ in the spectral sequence of the fibration $w^{\lambda} \mathbb{P}^n \rightarrow w^{\lambda} \mathbb{P}^n \times_{T^{\mu}} BT^{\mu} \rightarrow BT^{\mu}$.  The Leray-Hirsch theorem implies that $H^*_{T^{\mu}} (w^{\lambda} \mathbb{P}^n) \cong H^*(w^{\lambda} \mathbb{P}^n) \otimes H^*(BT^{\mu})$.  The inclusion $r: T^{\mu}_{\R} \rightarrow T^{\mu}$ induces an isomorphism $r^*: H^*_{T^{\mu}}(w^{\lambda} \mathbb{P}^n) \rightarrow H^*_{T^{\mu}_{\R}}(w^{\lambda} \mathbb{P}^n)$ and an isomorphism $H^*(BT^{\mu}) \cong H^*(BT^{\mu}_{\R})$.
\end{proof}

\begin{corollary} \label{inject to fixed points}
Fix $\mu$ and $\lambda$.  The inclusion of fixed points $\iota: (w^{\lambda} \mathbb{P}^n)^{T^{\mu}} \hookrightarrow w^{\lambda} \mathbb{P}^n$ induces an equivariant injection $\iota^*: H^*_{T^{\mu}}(w^{\lambda} \mathbb{P}^n) \hookrightarrow H^*_{T^{\mu}}((w^{\lambda} \mathbb{P}^n)^{T^{\mu}})$.  If $T^{\mu}_{\R}$ is the real torus then $\iota^*: H^*_{T^{\mu}_{\R}}(w^{\lambda} \mathbb{P}^n) \hookrightarrow H^*_{T^{\mu}_{\R}}((w^{\lambda} \mathbb{P}^n)^{T^{\mu}_{\R}})$ is also an injection.
\end{corollary}

\begin{proof}
Let $T^{\mu}_{\R}$ be the real subtorus of $T^{\mu}$.  The fixed points of $T^{\mu}_{\R}$ acting on $w^{\lambda} \mathbb{P}^n$ are the same as those of $T^{\mu}$.  The inclusion $r: T^{\mu}_{\R} \rightarrow T^{\mu}$ induces an isomorphism $r^*: H^*_{T^{\mu}}(w^{\lambda} \mathbb{P}^n) \rightarrow H^*_{T^{\mu}_{\R}}(w^{\lambda} \mathbb{P}^n)$.  Functoriality implies that $\iota^* \circ r^* = r^* \circ \iota^*$.  \cite[Corollary 2.2]{FP2} proves that $\iota^*$ is injective for the real torus $T^{\mu}_{\R}$ if $H^*_{T^{\mu}_{\R}} (w^{\lambda} \mathbb{P}^n)$ is a free $H^*(BT^{\mu}_{\R})$-module over $H^*(w^{\lambda} \mathbb{P}^n)$.
\end{proof}

For each one-dimensional $T^{\mu}$-orbit $O_{ij}$ the {\em weight} of the $T^{\mu}$-action on $O_{ij}$ is defined to be the Lie algebra of the stabilizer $\textup{Stab}(O_{ij})$, denoted $\mathfrak{t}_{ij}$.  If $[ce_j+de_i]$ is a point in $O_{ij}$ then
\[ (t_0,t_1, \ldots, t_n) \cdot [ce_j + de_i] = [c e_j + de_i]\]
if and only if there is $t \in \C^*$ such that as vectors
\[t_j^{\mu_j} ce_j + t_i^{\mu_i} de_i = t^{\lambda_j}ce_j + t^{\lambda_i}de_i.\]
Both $c$ and $d$ are nonzero so this is equivalent to the system of equations $t_j^{\mu_j} = t^{\lambda_j}$  and   $t_i^{\mu_i} = t^{\lambda_i}$.  Let $m_{ij}$ be $\textup{lcm} \{\lambda_i, \lambda_j \}$.  Then the weight on $O_{ij}$ is 
\[  \left(\frac{m_{ij}}{\lambda_i}\right) \mu_i t_i -   \left(\frac{m_{ij}}{\lambda_j} \right) \mu_j t_j .\]
(We use the same variables for $T^{\mu}$ as for $Lie(T^{\mu})$ and will assume $i<j$.)  

We state the GKM theorem as it applies to weighted projective space.  A more general description of GKM theory can be found in \cite{GKM} or the survey \cite{Tsurvey}.  

The cohomology ring $H^*_T(pt)$ is the symmetric algebra in the dual Lie algebra of $T$ and is typically written $\Z[t_0,t_1,\ldots,t_n]$.  GKM theory requires the topological space $X$ to have a finite number of $T$-fixed points.  This means $H^*_T(X^T) = \bigoplus_{x \in X^T} \Z[t_0,t_1,\ldots,t_n]$.  

The fixed points of weighted projective space are exactly the classes $[e_i]$ regardless of $\lambda$ or $\mu$.  
The class $p \in H^*_{T^{\mu}}((w^{\lambda} \mathbb{P}^n)^{T^{\mu}})$ is described using functional notation, meaning $p(x_i)$ denotes the $i^{th}$ polynomial of $p$.  The polynomial $p(x_i)$ is the {\em localization} of the class $p$ at the fixed point $x_i$.

\begin{theorem} \label{general GKM}
For each $\mu$ and $\lambda$ the rational equivariant cohomology of $w^{\lambda} \mathbb{P}^n$ is 
\[ H^*_{T^{\mu}}(w^{\lambda} \mathbb{P}^n; \Q) \cong \hspace{5in} \]
\[\hspace{.5in} \left\{p \in \bigoplus_{i=0}^n \Q[t_0,t_1,\ldots,t_n]: 
 \begin{array}{c} p(x_i) - p(x_j) \in \left\langle \left(\frac{m_{ij}}{\lambda_i}\right) \mu_i t_i -   \left(\frac{m_{ij}}{\lambda_j} \right) \mu_jt_j \right\rangle  \\
 \textup{ for each } 0 \leq i < j \leq n \end{array} \right\} .\]
\end{theorem}

The proof is immediate from \cite[Theorem 7.2]{GKM} since $w^{\lambda} \mathbb{P}^n$ is a $T^{\mu}$-CW space.  The ring structure of the localized classes is induced component-wise from the ring structure of the polynomial ring.  (Theorem \ref{general GKM} is also a corollary of the next theorem, using the map induced on equivariant cohomology with $\Q$ coefficients by the homomorphism that sends $(t_0,\ldots,t_n) \in T^{\mu}$ to $(t_0^{\mu_0},\ldots, t_n^{\mu_n}) \in T$.)

The modification for integral equivariant cohomology is subtle.  The main result in \cite{FP1} identifies the image of $H^*_{T^{\mu}}(w^{\lambda} \mathbb{P}^n)$ in the equivariant cohomology of the fixed points as long as the (real) isotropy subgroups of $T^{\mu}_{\R}$ are connected.  This is false when $\mu$ is not $(1,1,1,\ldots)$.  

\begin{theorem}
Let $\mu = (1,1,1,\ldots)$.  The equivariant cohomology ring of $w^{\lambda} \mathbb{P}^n$ is 
\[ H^*_{T}(w^{\lambda} \mathbb{P}^n; \Z) \cong \hspace{5in} \]
\[\hspace{.5in} \left\{p \in \bigoplus_{i=0}^n \Z[t_0,t_1,\ldots,t_n]: 
 \begin{array}{c} p(x_i) - p(x_j) \in \left\langle \left(\frac{m_{ij}}{\lambda_i}\right) t_i -   \left(\frac{m_{ij}}{\lambda_j} \right) t_j \right\rangle  \\
 \textup{ for each } 0 \leq i < j \leq n \end{array} \right\} .\]
\end{theorem}

\begin{proof}
We compute the equivariant cohomology with respect to the torus $T_{\R}$.  If $\mu = (1,1,1,\ldots)$ then each isotropy subgroup of $T$ is free in certain coordinates and is parametrized by the conditions that $t_i  = t^{\lambda_i}$ in the other coordinates (for $t \in \C^*$ free).  Hence each isotropy subgroup is connected.  With Proposition \ref{Tmu-CW}, this implies that \cite[Theorem 1.1]{FP1} holds: the image of $H^*_{T_{\R}}(w^{\lambda} \mathbb{P}^n) \stackrel{\iota^*}{\hookrightarrow} H^*_{T_{\R}}((w^{\lambda} \mathbb{P}^n)^{T_{\R}})$ agrees with the image of $H^*_{T_{\R}}(\overline{\bigcup O_{ij}}) \stackrel{\iota^*}{\hookrightarrow} H^*_{T_{\R}}((w^{\lambda} \mathbb{P}^n)^{T_{\R}})$.  

Each one-dimensional orbit $O_{ij}$ is $T_{\R}$-equivariantly contractible to $S^1$.  This implies $H^*_{T_{\R}}(O_{ij}) \cong \Z[\mathfrak{t}_{ij}^*]$.  The closure $\overline{O_{ij}}$ is a weighted projective line so $H^*_{T_{\R}}(\overline{O}_{ij})$ is a free $H^*_{T_{\R}}(pt)$-module over $H^*(\overline{O}_{ij})$ by Lemma \ref{equivariantly formal}.  Kawasaki proved the ordinary cohomology of $H^*(\overline{O}_{ij})$ is concentrated in even degree \cite{K}. Hence $H^*_{T_{\R}}(\overline{O}_{ij})$ is, too. The rest of the proof compares the Mayer-Vietoris sequence for the covering of $\overline{O_{ij}}$ by the standard open sets $O_{ij} \cup \{x_i\}$ and $O_{ij} \cup \{x_j\}$ to the long exact sequence for the pair $(\overline{O_{ij}}, \{x_i \cup x_j\})$.  It is identical to \cite[Theorem 7.2]{GKM}.   
\end{proof}

\section{A basis for $H^*_{T}(w^{\lambda} \mathbb{P}^n)$ and $H^*_{T^{\lambda}}(w^{\lambda} \mathbb{P}^n)$}

GKM theory comes equipped with natural combinatorial bases for equivariant cohomology called {\em canonical classes}.  In this section we describe explicitly the canonical classes of $H^*_{T^{\lambda}}(w^{\lambda} \mathbb{P}^n)$ when the weight of the torus action satisfies $\mu = \lambda$.   We will show that the canonical classes for weighted projective space are integer multiples of the canonical classes in $H^*_{T}(\mathbb{P}^n)$.  For ordinary projective space, these canonical classes are the geometric basis arising from the classical CW-decomposition.

Canonical classes can be defined in great generality \cite{GZ, T, GT}.  Not all varieties have canonical classes, but when canonical classes exist they generate equivariant cohomology \cite{GZ, T}.  We highlight the key properties for weighted projective space.
\begin{definition}
A canonical class $p_i^{\lambda}$ corresponding to the fixed point $x_i \in w^{\lambda} \mathbb{P}^n$ is an element $p_i^{\lambda} \in H^*_{T^{\mu}}(w^{\lambda} \mathbb{P}^n; \Q)$ that satisfies:
\begin{enumerate}
\item the localization $p_i^{\lambda}(x_i) = \prod_{j < i} \alpha_{ij}$ where $\alpha_{ij}$ is the weight on the orbit $O_{ij}$,
\item the class $p_i^{\lambda}$ is homogeneous of degree $i$, and 
\item the localization $p_i^{\lambda}(x_j) = 0$ if $j < i$.
\end{enumerate}
An {\em integral} canonical class $p_i^{\lambda} \in H^*_{T^{\mu}}(w^{\lambda} \mathbb{P}^n; \Z)$ satisfies $(2)$ and $(3)$ as well as
\begin{enumerate} \renewcommand{\labelenumi}{($\arabic{enumi}^\prime$)}
\item the localization $p_i^{\lambda}(x_i) = c \prod_{j < i} \alpha_{ij}$ for a constant $c \in \Z$. 
\end{enumerate}
\end{definition}
We denote the canonical classes of ordinary projective space by $p_i$.

\begin{lemma} \label{projective localizations}
Let $\mu = \lambda = (1,1,1, \ldots)$.  The class defined by
\[p_i(x_k) = \left\{ \begin{array}{ll} \prod_{j=0}^{i-1} (t_j - t_k) & \textup{ if } k \geq i, \textup{ and} \\  0 & \textup{ otherwise} \end{array} \right.\]
is an integral canonical class corresponding to $x_i$ in $H^*_T(\mathbb{P}^n)$.
\end{lemma}

\begin{proof}
The class $p_i$ satisfies all of the conditions of the definition of canonical classes by construction.  We need to confirm that it satisfies the GKM condition, namely that $p_i(x_k) - p_i(x_{k'}) \in  \langle t_{k'} - t_k \rangle$ for each $k' < k$.  This is implied directly by the definition if $k' < i$.  If $i \leq k'$ then $p_i(x_k)$ is exactly $p_i(x_{k'})$ with every instance of $t_{k'}$ replaced with $t_k$.  Hence the GKM condition holds.
\end{proof}

The canonical class $p_i$ is the geometric Schubert class corresponding to $x_i$ \cite{T}.  Rational canonical classes for arbitrary $\mu$ and $\lambda$ can be calculated in the same way but with more complicated coefficients.  Each integral canonical class is a multiple of a rational canonical class; the problem is to determine which multiple.

The next theorem explicitly constructs integral canonical classes for weighted projective spaces with $\mu = \lambda$.

\begin{theorem} \label{weighted basis}
Let $\jmath: \mathbb{P}^n \rightarrow w^{\lambda}\mathbb{P}^n$ be given by $\jmath[z_0,z_1,\ldots,z_n] = [z_0^{\lambda_0}, z_1^{\lambda_1}, \ldots, z_n^{\lambda_n}]$.  Then $\jmath$ induces an injective map
\[\jmath^*_T: H^*_{T^{\lambda}} (w^{\lambda}\mathbb{P}^n) \rightarrow H^*_T(\mathbb{P}^n).\]
Moreover there is a unique minimal canonical class $p_i^{\lambda} \in H^i_{T^{\lambda}} (w^{\lambda}\mathbb{P}^n)$ for each $i$, and $p_i^{\lambda}$ satisfies $\jmath^*(p_i^{\lambda}) = \kappa_i^{\lambda} p_i$ where $p_i$ is the ordinary cohomology canonical class and $\kappa_i^{\lambda}$ is the Kawasaki constant
\[\kappa_i^{\lambda} = \textup{lcm} \left\{ \frac{ \lambda_{j_0} \cdots \lambda_{j_i}}{\textup{gcd}\{ \lambda_{j_0}, \ldots, \lambda_{j_i} \}}: 0 \leq j_0 < j_1 < \cdots < j_i \leq n\right\}.\]
\end{theorem}

\begin{proof}
The map $\jmath$ is well-defined and equivariant with respect to the (ordinary) $T$-action on $\mathbb{P}^n$ and the (weighted) $T^{\lambda}$-action on $w^{\lambda} \mathbb{P}^n$, in the sense that 
\[\jmath((t_0,\ldots,t_n) \cdot [z_0,\ldots, z_n]) = (t_0,\ldots,t_n)^{\lambda} \cdot \jmath([z_0,\ldots,z_n]).\]
Hence $\jmath$ induces a map on equivariant cohomology as in \cite[Lecture 2]{F}.
The set of $T^{\lambda}$-fixed points in $w^{\lambda}\mathbb{P}^n$ is the same as the set of $T$-fixed points in $\mathbb{P}^n$.  The inclusion $\iota$ of fixed points into weighted and ordinary projective space commutes with $\jmath$.  This gives a commutative diagram
\[\begin{array}{ccc}
H^*_{T^{\lambda}}(w^{\lambda} \mathbb{P}^n) & \stackrel{\jmath^*_T}{\longrightarrow} & H^*_T(\mathbb{P}^n)  \\
\downarrow \iota^* & & \downarrow \iota^*  \\
H^*_{T^{\lambda}}\left( (w^{\lambda} \mathbb{P}^n)^{T^{\lambda}} \right) & \stackrel{\jmath^*_T}{\longrightarrow} & H^*_{T}\left(( \mathbb{P}^n)^{T}\right). \end{array} \]
\cite[Lecture 3, Example 2.2]{F} shows that the map on the bottom row restricts the algebra homomorphism 
\[\jmath^*_T: \bigoplus_{i=0}^n \Z[t_0,\ldots,t_n] \rightarrow \bigoplus_{i=0}^n \Z[t_0,\ldots,t_n]\] 
given on the generators by $\jmath^*_T(t_i) = \lambda_i t_i$.  Hence $\jmath^*_T$ is injective on $H^*_{T^{\lambda}}\left( (w^{\lambda} \mathbb{P}^n)^{T^{\lambda}} \right)$.  Corollary \ref{inject to fixed points} proved $\iota^*$ is injective, so $\jmath^*_T$ is injective in the top row as well.

Lemma \ref{equivariantly formal} showed that weighted projective space is equivariantly formal for all $\lambda$ and all $T^{\mu}$-actions.  Hence we have surjections $H^*_{T^{\lambda}}(w^{\lambda} \mathbb{P}^n) \rightarrow H^*(w^{\lambda} \mathbb{P}^n)$ and $H^*_{T}(\mathbb{P}^n) \rightarrow H^*(\mathbb{P}^n)$.  The equivariant map $\jmath^*_T$ is a $\Z[t_0,\ldots,t_n]$-algebra homomorphism so the following diagram commutes
\[\begin{array}{ccc}
H^*_{T^{\lambda}}(w^{\lambda} \mathbb{P}^n) & \stackrel{\jmath^*_T}{\longrightarrow} & H^*_T(\mathbb{P}^n)  \\
\downarrow  & & \downarrow  \\
H^* \left( w^{\lambda} \mathbb{P}^n \right) & \stackrel{\jmath^*}{\longrightarrow} & H^*\left( \mathbb{P}^n\right). \end{array} \]
The induced map in ordinary cohomology $\jmath^*$ is the Kawasaki homomorphism from \cite{K}. Kawasaki showed that if $p_i^{\lambda}$ denotes the generator of $H^i( w^{\lambda} \mathbb{P}^n )$ (respectively $p_i$ and $H^*(\mathbb{P}^n)$) then $\jmath^*(p_i^{\lambda}) = \kappa_i^{\lambda} p_i$.

The $\Z[t_0,\ldots,t_n]$-module $H^*_{T^{\lambda}}(w^{\lambda} \mathbb{P}^n)$ is naturally isomorphic by equivariant formality to $H^*(w^{\lambda} \mathbb{P}^n) \otimes \Z[t_0,\ldots,t_n]$.  Let $p_i^{\lambda}$ also denote the element of $H^*_{T^{\lambda}}(w^{\lambda} \mathbb{P}^n)$ that corresponds under this isomorphism to $p_i^{\lambda} \otimes 1$.  Similarly use $p_i$ to denote the canonical class in $H^*_T(\mathbb{P}^n)$, which corresponds to $p_i \otimes 1$.  The commutativity of the diagram {\em as module homomorphisms} allows us to conclude $\jmath^*_T(p_i^{\lambda}) = \kappa_i^{\lambda} p_i$ in equivariant cohomology too.  The injectivity of $\jmath^*_T$ tells us this is the unique equivariant class in the preimage of $\kappa_i^{\lambda} p_i$.

Recall that $\jmath^*_T \circ \iota^* = \iota^* \circ \jmath^*_T$ and that $\kappa_i^{\lambda} p_i$ is a canonical class in $H^*_T(\mathbb{P}^n)$.  Hence the localization $p_i^{\lambda}(x_k) = \kappa_i^{\lambda} \prod_{j<i} \left( \frac{t_j}{\lambda_j} -\frac{t_k}{\lambda_k} \right) $ in $H^*_{T^{\lambda}}(pt)$.  To show $p_i^{\lambda}$ is a canonical class in $H^*_{T^{\mu}}(w^{\lambda} \mathbb{P}^n)$ we prove $\textup{lcm}\{\lambda_j, \lambda_k\}$ divides $\kappa_i^{\lambda}$ for all possible $i,j,k$.  Our strategy is to show that for each $k$ and each $i$ the integer $\lambda_k$ divides at least one term $\frac{ \lambda_{j_0} \cdots \lambda_{j_i}}{\textup{gcd}\{ \lambda_{j_0}, \ldots, \lambda_{j_i} \}}$.  Note that $\lambda_k$ divides $\frac{\lambda_k \lambda_j}{\textup{gcd}\{\lambda_k,\lambda_j\}}$ for each $j$ since $\frac{\lambda_j}{\textup{gcd}\{\lambda_k,\lambda_j\}}$ is a positive integer.  Similarly $\textup{gcd}\{ \lambda_{j_0}, \ldots, \lambda_{j_i} \}$ divides $\lambda_{j_i}$ for any choice of $\lambda_{j_i}$.  If $\lambda_k$ divides $\frac{ \lambda_{j_0} \cdots \lambda_{j_{i-1}}}{\textup{gcd}\{ \lambda_{j_0}, \ldots, \lambda_{j_{i-1}} \}}$ then  we conclude $\lambda_k$ divides $\frac{ \lambda_{j_0} \cdots \lambda_{j_i}}{\textup{gcd}\{ \lambda_{j_0}, \ldots, \lambda_{j_i} \}}$.  By induction $\lambda_k$ divides $\kappa_i^{\lambda}$ for each $i$ and $k$.  Hence $\textup{lcm}\{\lambda_j, \lambda_k\}$ divides $\kappa_i^{\lambda}$ for all $i,j,k$.

Thus $p_i^{\lambda}$ is a canonical class for $x_i$.  If there is another canonical class for $x_i$, say $q_i$, then the image of $q_i$ in ordinary cohomology is sent to $c (\kappa_i^{\lambda} p_i)$ for some integer $c$.  Hence $p_i^{\lambda}$ is the minimal canonical class for $x_i$.
\end{proof}

\section{The rings $H^*_{T}(\mathbb{P}^n)$ and $H^*_{T^{\lambda}}(w^{\lambda} \mathbb{P}^n)$}

In this section we obtain an explicit formula for the structure constants in $H^*_{T^{\lambda}}(w^{\lambda}\mathbb{P}^n)$.  We begin with ordinary projective space, which by the previous section is the key to the weighted case.  Amazingly for an object as long-studied as $\mathbb{P}^n$, this result appears to be new.

The divided difference operator $\partial_j$ is the rational operator on the polynomial ring $\Z[t_1,\ldots, t_n]$ defined by
\[ \partial_j p(t_1, \ldots, t_n) = \frac{p(\ldots, t_{j+1}, t_j, \ldots) - p(\ldots, t_j, t_{j+1}, \ldots)}{t_j-t_{j+1}}.\]
(Some authors use different sign conventions.)

\begin{theorem} \label{ordinary equivariant projective}
Suppose $0 \leq i \leq j < k \leq i+j, n$.  The nonzero structure constants of $H^*_T(\mathbb{P}^n)$ are $c_{ij}^j = p_i(x_j)$ and
\[c_{ij}^k = \partial_{k-1} \partial_{k-2} \cdots \partial_{j+2} \partial_{j+1} \partial_j p_i(x_j).\]
\end{theorem}

\begin{proof}
The proof is by induction on $k$ with separate cases for $k<j$ and $k>j$.  We begin by writing the product $p_i p_j = \sum_{k=0}^{n} c_{ij}^k p_k$ in terms of the basis of canonical classes, for some coefficients $c_{ij}^k \in \C[t_0,t_1,\ldots, t_n]$.  If $k > i+j$ then the coefficient $c_{ij}^k=0$ since $p_ip_j$ is homogeneous of degree $i+j$ and each $p_k$ has degree $k$.  

The canonical class $p_k$ vanishes at the localization $p_k(x_j)=0$ when $k > j$.  Hence localizing the product $p_ip_j$ at $x_k$ when $k < j$ gives
\[p_i(x_k)p_j(x_k) = 0.\]
If $0<j$ then localizing at $x_0$ gives $\sum_{k=0}^{i+j} c_{ij}^k p_k(x_0) = c_{ij}^0 p_0(x_0) = 0$.  The definition of canonical classes states that $p_0(x_0) \neq 0$ so we conclude that $c_{ij}^0 = 0$.  Assume as the inductive hypothesis that $c_{ij}^0, c_{ij}^1, \ldots, c_{ij}^{k-1}$ are all zero and that $k<j$.  Then localizing at $x_k$ gives $p_i(x_k)p_j(x_k) = 0 = c_{ij}^k p_k(x_k)$.  Again we conclude $c_{ij}^k=0$.  By induction $c_{ij}^k=0$ for all $k<j$.

Now consider the cases $k \geq j$.   Localizing at $x_j$ gives
\[p_i(x_j)p_j(x_j) = c_{ij}^j p_j(x_j).\]
This determines $c_{ij}^j$.  We compute $c_{ij}^{j+1}$ in the same way.  First we localize at $x_{j+1}$:
\[p_i(x_{j+1})p_j(x_{j+1}) = p_i(x_j)p_j(x_{j+1}) + c_{ij}^{j+1} p_{j+1}(x_{j+1}).\]
Next we solve for $c_{ij}^{j+1}$:
\[c_{ij}^{j+1} = \frac{p_i(x_{j+1})p_j(x_{j+1}) -  p_i(x_j)p_j(x_{j+1})}{p_{j+1}(x_{j+1})}.\]
Lemma \ref{projective localizations} shows that $\frac{p_{j+1}(x_{j+1})}{p_j(x_{j+1})} = t_j - t_{j+1}$.  Hence we may simplify
\[c_{ij}^{j+1} = \frac{p_i(x_{j+1}) -  p_i(x_j)}{t_j - t_{j+1}}.\]
Lemma \ref{projective localizations} also shows that $p_i(x_{j+1})$ is obtained from $p_i(x_j)$ by exchanging the roles of $t_j$ and $t_{j+1}$, so $c_{ij}^{j+1} = \partial_j p_i(x_j)$.

Assume as the inductive hypothesis that the formula holds through $k$.  If $k \geq j+1$ then localizing at $x_k$ gives
\[p_i(x_k) p_j(x_k) = p_i(x_j)p_j(x_k) + \cdots + c_{ij}^k p_k(x_k).\]
This can be rearranged as
\[c_{ij}^k = \frac{p_i(x_k) p_j(x_k) - p_i(x_j)p_j(x_k) - \cdots -  c_{ij}^{k-1} p_{k-1}(x_k)}{p_k(x_k)}.\]
Similarly, we may write
\begin{equation}\label{recursion}
c_{ij}^{k+1} =  \frac{p_i(x_{k+1}) p_j(x_{k+1}) - p_i(x_j)p_j(x_{k+1}) - \cdots -  c_{ij}^{k-1} p_{k-1}(x_{k+1})}{p_{k+1}(x_{k+1})} - \frac{c_{ij}^k p_k(x_{k+1})}{p_{k+1}(x_{k+1})}.\end{equation}
The denominator factors as $p_{k+1}(x_{k+1}) = p_k(x_{k+1}) (t_k - t_{k+1})$ by Lemma \ref{projective localizations}.  The second term in the sum in Equation \eqref{recursion} is thus $c_{ij}^k/(t_k-t_{k+1})$. 

If $p$ is a polynomial, let $s_k p$ denote the polynomial obtained by exchanging the variables $t_k$ and $t_{k+1}$ in $p$.  

We claim the first term in Equation \eqref{recursion} is $\frac{s_k c_{ij}^{k}}{t_k-t_{k+1}}$.  The variables $t_k, t_{k+1}$ do not appear in $p_i(x_j)$ by Lemma \ref{projective localizations}.  The divided difference operators $\partial_j, \partial_{j+1}, \ldots, \partial_{k-2}$ do not introduce the variables $t_k, t_{k+1}$, so by the inductive hypothesis neither $t_k$ nor $t_{k+1}$ appears in the coefficients $c_{ij}^j, c_{ij}^{j+1}, \ldots, c_{ij}^{k-1}$.  For each $k' \leq k$ the polynomial $p_{k'}(x_k)$ is $s_k p_{k'}(x_{k+1})$ by Lemma \ref{projective localizations}.  So the numerator of the first term in Equation \eqref{recursion} is the numerator of $c_{ij}^k$ with the roles of $t_k$ and $t_{k+1}$ exchanged.  Since the denominator is $p_{k+1}(x_{k+1}) = p_k(x_{k+1}) (t_k - t_{k+1})$ we may write
\[c_{ij}^{k+1} = \frac{s_k c_{ij}^k - c_{ij}^k}{t_k - t_{k+1}} = \partial_k c_{ij}^k.\]
The claim follows by induction.
\end{proof}

We remark that the coefficients $c_{ij}^k$ are not in general Schubert polynomials, neither as defined by Lascoux-Schutzenberger \cite{LS} nor as defined by Bernstein-Gelfand-Gelfand \cite{BGG}.

The next corollary gives a direct and concise combinatorial proof that the polynomials $c_{ij}^k$ have nonnegative coefficients when written in terms of the variables $\alpha_1 = t_1 - t_2, \alpha_2 = t_2 - t_3, \ldots, \alpha_{n-1} = t_{n-1} - t_n$.  This extends work of Graham, who proved the analogous statement for the equivariant cohomology of generalized flag varieties with respect to the basis of Schubert classes \cite{G}.  (Positivity for $H^*_T(\mathbb{P}^n)$ can also be proven indirectly by realizing the equivariant cohomology of projective space as a subring of the equivariant cohomology of the flag variety and then applying Graham's result.)

\begin{corollary}
The coefficients $c_{ij}^k$ are polynomials with nonnegative coefficients when written in the variables $\alpha_1 = t_1 - t_2, \alpha_2 = t_2 - t_3, \ldots, \alpha_{n-1} = t_{n-1} - t_n$.
\end{corollary}

\begin{proof}
The transpositions $s_j$ act on the variables $\alpha_1, \alpha_2, \ldots$ by
\[s_j (\alpha_k) = \left\{ \begin{array}{ll} \alpha_k+\alpha_j & \textup{ if } j = k \pm 1, \\
	- \alpha_k & \textup{ if } j = k, \textup{ and} \\
	\alpha_k & \textup{ otherwise.} \end{array} \right.\]
The polynomial $p_i(x_j)$ is written in terms of the variables $\alpha_1, \alpha_2, \ldots$ as 
\[p_i(x_j) =(\alpha_1 + \alpha_{2} + \cdots + \alpha_{j-1})(\alpha_2 +\alpha_3+ \cdots + \alpha_{j-1}) \cdots (\alpha_{i-1} + \alpha_i + \cdots + \alpha_{j-1}).\]  
Thus $p_i(x_j)$ is a polynomial in $\alpha_{j-1}$ with coefficients in $\mathbb{Z}_{\geq 0}[\alpha_1, \ldots, \alpha_{j-2}]$.

The proof proceeds by induction on $k$.  The inductive hypothesis is that $c_{ij}^{k}$ can be written as a polynomial in $\alpha_{k-1}$ whose coefficients are in $\mathbb{Z}_{\geq 0}[\alpha_1, \ldots, \alpha_{k-2}]$.  Denote this polynomial $q_k(\alpha_{k-1})$.  Then $s_k$ fixes each coefficient of $q_k$, so 
\[ \partial_k \left(c_{ij}^k \right) = \frac{s_k (q_k(\alpha_{k-1})) - q_k(\alpha_{k-1})}{\alpha_k} = \frac{q_k (\alpha_{k-1} + \alpha_k) - q_k(\alpha_{k-1})}{\alpha_k}.\]
Collecting the coefficients of powers of $\alpha_k$ lets us write 
\[q_k (\alpha_{k-1} + \alpha_k) = q_k(\alpha_{k-1}) + \alpha_k q_{k+1}(\alpha_k)\] 
for some polynomial $q_{k+1}(\alpha_k)$ with coefficients in  $\mathbb{Z}_{\geq 0}[\alpha_1, \alpha_2, \ldots, \alpha_{k-1}]$.  Since $\partial_k c_{ij}^k = c_{ij}^{k+1}$ it follows by induction that all $c_{ij}^k$ are nonnegative polynomials in the variables $\alpha_1, \alpha_2,\ldots$.
\end{proof}

Finally we compute the equivariant structure constants of $w^{\lambda} \mathbb{P}^n$.

\begin{corollary}
Let $c_{ij}^k$ denote the structure constants in $H^*_T(\mathbb{P}^n)$.  Let $\kappa_i^{\lambda}$ denote the integers from Theorem \ref{weighted basis}, so that each minimal canonical class $p_i^{\lambda} \in H^*_{T^{\lambda}}(w^{\lambda} \mathbb{P}^n)$ satisfies $p_i^{\lambda} = \kappa_i^{\lambda} p_i$.  The structure constants $c_{ij}^{k,\lambda}$ in the product $p_i^{\lambda} p_j^{\lambda} = \sum c_{ij}^{k,\lambda} p_k^{\lambda}$ are given by
\[c_{ij}^{k,\lambda} = \frac{\kappa_i^{\lambda}\kappa_j^{\lambda}}{\kappa_k^{\lambda}}c_{ij}^k.\]
\end{corollary}

\begin{proof}
If $p_i p_j = \sum c_{ij}^k p_k$ in $H^*_T(\mathbb{P}^n)$ then 
\[(\kappa_i^{\lambda}p_i)(\kappa_j^{\lambda}p_j) = \sum \left( \frac{\kappa_i^{\lambda}\kappa_j^{\lambda}}{\kappa_k^{\lambda}}c_{ij}^k \right) (\kappa_k^{\lambda}p_k)\]
inside the ring $ \bigoplus_{i=0}^n \Z[t_0,t_1,\ldots,t_n]$.  
\end{proof}

\end{document}